\documentclass[12pt,a4paper]{amsart}

\usepackage {amsfonts}
\usepackage[english]{babel}
\def\b{\beta}
\def\g{\gamma}
\def\G{\Gamma}
\def\d{\delta}
\def\a{\alpha}
\def\p{\varphi}
\def\e{\varepsilon}
\def\l{\lambda}
\def\L{\Lambda}
\def\t{\theta}
\def\la{\langle}
\def\ra{\rangle}
\def\R{{\mathbb R}}
\def\C{{\mathbb C}}
\def\N{{\mathbb N}}
\def\Z{{\mathbb Z}}
\def\cM{{\mathcal M}}
\def\hr{\frak r}

\def\hR{\frak R}

\def\o{\omega}
\def\p{\varphi}
\def\bs{~\hfill\rule{7pt}{7pt}}
\def\F{Fourier quasicrystal}

\DeclareMathOperator{\dist}{dist}
\DeclareMathOperator{\supp}{supp}
\DeclareMathOperator{\Lin}{Lin}
\DeclareMathOperator{\di}{dim}

\newtheorem{Th}{Theorem}
\newtheorem*{Cor}{Corollary}
\newtheorem{Lem}{Lemma}

\begin{document}

\title{Large Fourier Quasicrystals and Wiener's Theorem}

\author{S.Yu.Favorov}

\address{Sergii Favorov,
\newline\hphantom{iii}  Karazin's Kharkiv National University
\newline\hphantom{iii} Svobody sq., 4,
\newline\hphantom{iii} 61022, Kharkiv, Ukraine}
\email{sfavorov@gmail.com}

\maketitle {\small
\begin{quote}
\noindent{\bf Abstract.}
We find new simple conditions for support of a discrete measure on  Euclidean space to be a finite union of translated lattices. The arguments are based on a local analog of Wiener's Theorem on absolutely convergent trigonometric series and theory of almost periodic functions.
\medskip


\end{quote}
}

\section{introduction}

Denote by $S(\R^d)$ the Schwartz space of test functions $\p\in C^\infty(\R^d)$ with finite norms
 $$
  p_m(\p)=\sup_{\R^d}(\max\{1,|x|\})^m\max_{|k_1|+\dots+|k_d|\le m} |D^k(\p(x))|,\quad m=0,1,2,\dots,
 $$
  $k=(k_1,\dots,k_d)\in(\N\cup\{0\})^d,\  D^k=\partial^{k_1}_{x_1}\dots\partial^{k_d}_{x_d}$. These norms generate a topology on $S(\R^d)$, and elements of the space $S'(\R^d)$ of continuous linear functionals on $S(\R^d)$ are called tempered distributions.
  For every tempered distribution $f$ there exist $C>0$ and $m\in\N\cup\{0\}$ such that for all $\p\in S(\R^d)$
$$
                           |f(\p)|\le Cp_m(\p).
$$
Moreover, this estimate is sufficient for distribution $f$ to be in $S'(\R^d)$
(see \cite{V}, Ch.3).

The Fourier transform of a tempered distribution $f$ is defined by the equality
$$
\hat f(\p)=f(\hat\p)\quad\mbox{for all}\quad\p\in S(\R^d),
$$
where
$$
   \hat\p(y)=\int_{\R^d}\p(x)\exp\{-2\pi i\la x,y\ra\}dm_d(x)
 $$
is the Fourier transform of the function $\p$. Note that the Fourier transform of every tempered distribution is also a tempered distribution. Here we consider only the case when $f$ and $\hat f$ are measures on $\R^d$.
\medskip

To formulate  the results of the paper we need some notions and definitions.

Set $B(x,r)=\{y\in\R^d:\,|y-x|<r\}$, $B(r)=B(0,r)$. We say that $E\subset\R^d$ is {\it a full-rank lattice} if  $E=T\Z^d$ for some nondegenerate linear operator $T$ on $R^d$; $E$ is {\it relatively dense} if there exists $R<\infty$ such that $B(x,R)\cap E\neq\emptyset$ for each $x\in\R^d$; $E$ is {\it discrete}  if $E\cap B(x,1)$ is finite for all $x\in\R^d$;  $E$ is {\it uniformly discrete} if $|x-x'|\ge\e>0$ for all $x,\,x'\in E, x\neq x'$.

 A Radon measure $\mu$ is discrete (uniformly discrete) if its support is discrete (uniformly discrete), $\mu$ is {\it translation bounded} if its variations $|\mu|$  are uniformly bounded on balls of radius 1, and $\mu$ is  {\it slowly increasing} if $|\mu|(B(r))$ grows at most polynomially as $r\to\infty$. Note that every translation bounded measure is slowly increasing, and every slowly increasing measure belongs to $S'(\R^d)$.

Following \cite{LO2} we say that $\mu$ is a {\it\F} if both measures $\mu$, $\hat\mu$ are slowly increasing atomic measures.  More precisely, we will suppose that
$$
\mu=\sum_{\l\in\L}a_\l\d_\l, \quad\hat\mu=\sum_{\g\in\G}b_\g\d_\g, \quad \log\left[\sum_{|\l|<r}|a_\l|+\sum_{|\g|<r}|b_\g|\right]= O(\log r),\ r\to\infty,
$$
where $\d_z$ means the unit mass at the point $z$ and the sets $\L,\,\G$ are countable. In the paper we usually consider the case
of discrete {\it support} $\L$ and countable {\it spectrum} $\G$ of  $\mu$. We say that $\mu$ is a {\it large \F} if, in addition, $\inf_{\l\in\L}|a_\l|>0$.
The simplest example of a \F\ is the measure $\mu_0=\sum_{k\in\Z^d}\d_k$. By the Poisson formula we have $\hat\mu_0=\mu_0$.

Such measures are the main object in the theory of \F s (see \cite{C}-\cite{Mo}). The corresponding notions were inspired by experimental discovery  of non-periodic atomic structures with diffraction patterns consisting of spots, which was made in the mid '80s.

\medskip
In the paper we are interested in the cases when support $\L$ is a subset of a finite union of translated full-rank lattices. Several very interesting results of this type were obtained by N.Lev, A.Olevskii in \cite{LO1}, \cite{LO2}. Here it is one of them.
 \begin{Th}[\cite{LO2}]
Let $\mu$ be a \F\ on $\R^d$ with discrete spectrum $\G$ and uniformly discrete set of differences $\L-\L$. Then  $\L$ is a subset of a finite union of translates of a single full-rank lattice $L$, and $\G$ is a subset of a finite union of translates of the conjugate lattice.

Also, there exists a \F\  with countable spectrum $\G$ such that $\L-\L$ is uniformly discrete, but $\L$ is not contained in a finite union of translates of any lattice.
\end{Th}

We prove the following theorem,  which amplifies the previous one

 \begin{Th}\label{3}
Let $\mu$ be a large \F\ on $\R^d$ with a discrete set of differences $\L-\L$. Then  $\L$ is a  finite union of translates of a single full-rank lattice $L$.
\end{Th}
In the next section we will give a generalization of Theorem \ref{3} for pairs of \F s.
\medskip

Results of another type of results were obtained by Y.Meyer and A.Cordoba.
\begin{Th}[Y.Meyer \cite{M1}]\label{1}
Let $\mu=\sum_{\l\in\L} a_\l\d_\l$ be a  measure with discrete support $\L$ and $a_\l\in S$ for some finite set $S\subset\C\setminus\{0\}$.
If $\mu\in S'(\R)$ and if its Fourier transform $\hat\mu$ is a translation bounded measure on $\R$, then $\L=E\triangle\bigcup_{j=1}^N(\a_j\Z+\b_j)$, where $\a_j>0,\ \b_j\in\R$, and the set $E$ is finite.
\end{Th}
Here $A\triangle B$ means the symmetrical difference between $A$ and $B$.
\medskip

 In \cite{K} M.Kolountzakis extended the above theorem to measures on $\R^d$. He replaced the condition "the measure $\hat\mu$ is translation bounded" with the weaker one
 \begin{equation}\label{m}
|\hat\mu|(B(r))=O(r^d)\quad \mbox{as} \quad r\to\infty.
\end{equation}
He also found a condition for support $\mu$ to be a finite union of several full-rank lattices. His result is very close to Cordoba's  one:
  \begin{Th}[\cite{C}]\label{2}
 Let  $\mu$ be a uniformly discrete \F\ on $\R^d$ with $a_\l$  belonging to a finite set $F$.  If the measure $\hat\mu$ is translation bounded,  then $\L$  is a finite union of translates of several, possibly incommensurable, full-rank lattices.
\end{Th}

In paper \cite{F2} we replaced the conditions "$a_\l$ from a finite set" by "$|a_\l|$ from a finite set" and the condition "$\hat\mu$ is translation bounded" by (\ref{m}).

Now we obtain the following theorem.
\begin{Th}\label{4}
 Let $\mu$  be a uniformly discrete large \F\ on $\R^d$,  and let $\hat\mu$  satisfy (\ref{m}). Then  $\L$  is a finite union of translates of several disjoint full-rank lattices.
 \end{Th}
  The proof is based on an analog of Wiener's Theorem on Fourier series:
\begin{Th}[N.Wiener, see, for example, \cite{Z}, Ch.VI]
Let $F(t)=\sum_{n\in\Z}c_ne^{2\pi int}$  be an absolutely convergent Fourier series,  and $h(z)$  be a holomorphic function on a neighborhood of  the set $\overline{\{F(t):\,t\in [0,1]\}}$. Then the function $h(F(t))$ admits  an absolutely convergent Fourier series expansion as well.
 \end{Th}
 Denote by $W$ the class of absolutely convergent series $\sum_n c_ne^{2\pi i\la x,\g_n\ra},\ \g_n\in\R^d$. We prove the following theorem
 \begin{Th}\label{5}
 Let  $K\subset\C$ be an arbitrary compact,  $h(z)$ be a holomorphic function on a neighborhood of $K$, and $f\in W$. Then there is a function $g\in W$ such that if $f(x)\in K$ then $g(x)=h(f(x))$.
 \end{Th}
 If $K=\overline{f(\R^d)}$ then we obtain the global Wiener's Theorem for almost periodic functions.
 \medskip

 Finally, we prove an analog of Wiener's Theorem for \F s:
 \begin{Th}\label{7}
 Let $\mu=\sum_{\l\in\L}a_\l\d_\l$  be a uniformly discrete \F\ on $\R^d$ with the translation bounded measure $\hat\mu$.  Suppose that $h(z)$ is a holomorphic function on a neighborhood of the closure of the set ${\{a_\l:\,\l\in\L\}}$. Then  $\nu=\sum_{\l\in\L}a_\l h(a_\l)\d_\l$  is also \F\ with  translation bounded $\hat\nu$.
 \end{Th}
 \begin{Cor}
  Let $\mu=\sum_{\l\in\L}a_\l\d_\l$  be a uniformly discrete \F\ on $\R^d$ with translation bounded measure $\hat\mu$. Let $A$ be a connected component of the set  $\overline{{\{a_\l:\,\l\in\L\}}}$ such that $0\notin A$. Then  $\{\l\in\L:\,a_\l\in A\}$  is a finite union of translates of several full-rank lattices.
 \end{Cor}

\section{Large quasicrystals and almost periodic measures}

\bigskip

We recall some definitions  related to the notion of almost periodicity (see, for example, \cite{M}).

A continuous function $f$ on $\R^d$ is {\it almost periodic}, if
for every  $\e>0$ the set of $\e$-almost periods of $f$
  $$
  \{\tau\in\R^d:\,\sup_{x\in\R^d}|f(x+\tau)-f(x)|<\e\}
  $$
  is a relatively dense set in $\R^d$.

  For every almost periodic function $f$ there exists the limit uniform in $x\in\R^d$
$$
\cM(f)=\lim_{r\to\infty}\frac{1}{\o_d r^d}\int_{B(x,r)}f(y)dm_d(y),
$$
where $\o_d$ is volume of the unit ball. The Fourier coefficient of $f$ at frequency $\g\in\R^d$ is
$$
a_f(\g)=\cM\left(f(x)e^{-2\pi i\la x,\g\ra}\right).
$$
 For every $f$ only countably many $a_f(\g)$ do not vanish,
and the set $\{a_f(\g)\}$ defines the function $f$ uniquely.

 It can be easily proved that every function $f=\sum_n c_ne^{2\pi i\la x,\g_n\ra}\in W$ is almost periodic. We set $\|f\|_W=\sum_n|c_n|$.

\noindent Since $\cM\left(e^{2\pi i\la x,\g\ra}e^{2\pi i\la x,-\g'\ra}\right)=\begin{cases}1, & \g=\g', \\0, & \g\neq\g',\end{cases}$\qquad
we get $a_f(\g)=\begin{cases}c_n, & \g=\g_n, \\0, & \g\neq\g_n.\end{cases}$

\noindent Moreover,  $\hat f=\sum_n c_n\d_{\g_n}$. Indeed, if $\psi\in S(\R^d)$, then
$$
   (\hat f,\psi)=(f,\hat\psi)=\sum_nc_n\int_{\R^d}e^{2\pi i\la y,\g_n\ra}\hat\psi(y)dy=\sum_n c_n\psi(\g_n).
$$
Next, a (complex) measure $\mu$ on $\R^d$ is {\it almost periodic}, if
for every continuous function $\psi$ on $\R^d$ with compact support
 the function $(\psi\star\mu)(x)$ is almost periodic in $x\in\R^d$. It is clear that  support of every almost periodic measure is relatively dense in $\R^d$.

A discrete set $\L$ is almost periodic, if measure $\sum_{\l\in\L}\d_\l$ is almost periodic.

\begin{Th}[L.Ronkin \cite{R}]\label{R}
Every almost periodic measure is translation bounded.
\end{Th}
A connection between almost periodicity of measure and properties of its Fourier transform was obtained by Y.Meyer:
\begin{Th} [\cite{M}]
Let $\mu$ and its Fourier transform $\hat\mu$ be translation bounded measures.
Then $\mu$ is almost periodic iff spectrum of $\mu$  is countable.
\end{Th}
Here we need an extension of this result.
\begin{Th}\label{9}
 Let $\mu$ be a measure from $S'(\R^d)$, and let $\hat\mu$ be a slowly increasing measure. Then $\mu$ is almost periodic  iff it is translation bounded and has a countable spectrum.
\end{Th}
First we prove some lemmas.
\begin{Lem}\label{L1}
 Let $\mu$ be a measure from $S'(\R^d)$ with countable spectrum, let $\hat\mu$ be a slowly increasing measure, and let $\psi\in S(\R^d)$.
 Then convolution $\psi\star\mu$ belongs to $W$.
\end{Lem}
{\bf Proof}. Let $\hat\mu=\sum_{\g\in\G}b_\g\d_\g$. We have
\begin{equation}\label{a}
(\psi\star\mu)(s)=\int\psi(s-x)d\mu(x)=\int\hat\psi(y)e^{2\pi i\la y,s \ra}d\hat\mu(y)=\sum_{\g\in\G}b(\g)\hat\psi(\g)e^{2\pi i\la s,\g \ra}.
\end{equation}
Since $\hat\psi\in S(\R^d)$, we get
 $$
\int|\hat\psi(y)|d|\hat\mu|(y) \le\int_0^\infty C\min\{1,\,t^{-N}\}dM(t) =C\lim_{t\to\infty}\frac{M(t)}{t^N}+CN\int_1^\infty\frac{M(t)}{t^{N+1}}dt,
$$
where $M(t)=|\hat\mu|(B(t))$. For appropriate $N$ the latter integral is finite. Therefore the sum in (\ref{a}) converges absolutely. \bs
\begin{Lem}\label{L2}
 Let $\mu$ be a translation bounded measure, and let $\hat\mu$ be a slowly increasing measure.
 Then for each $\l\in\R^d$ the limit
 $$
    \lim_{R\to\infty}\frac{1}{\o_d R^d}\int_{B(R)}e^{-2\pi i\la x,\l\ra}d\mu(x)
 $$
exists and equals $\hat\mu(\{\l\})$.
\end{Lem}
{\bf Proof}. Pick an arbitrary  $\e>0$. Let $\psi(|x|)$ be $C^\infty$-differentiable function on $\R^d$ such that $\psi(t)=1/\o_d$ for $t<1-\e$, $\psi(t)=0$ for $t>1+\e$, $0\le\psi(t)\le1/\o_d$ for $1-\e \le t\le 1+\e$, and $\int_{\R^d}\psi(|x|)dx=1$. We have
$$
   \frac{1}{R^d}\int_{\R^d}\psi(|x|/R)e^{-2\pi i\la x,\l\ra}d\mu(x)= \{R^{-d}\psi(|x/R|)\mu(x)\}\hat\, (\l)=\hat\psi(R(\l-y))\star\hat\mu(y),
 $$
 then
\begin{equation}\label{j}
\hat\psi(R(\l-y))\star\hat\mu(y)=\int_{|\l-y|<1}\hat\psi(R(\l-y))d\hat\mu(y)
+\int_{|\l-y|\ge1}\hat\psi(R(\l-y))d\hat\mu(y).
 \end{equation}
  Note that $\hat\psi(0)=1$ and $\hat\psi(R(\l-y))\to0$ for $\l\neq y$ as $R\to\infty$. By the Dominated Convergence Theorem, the first integral in (\ref{j}) tends to $0$ if $\hat\mu(\{\l\})=0$ and to $\hat\mu(\{\l\})$ otherwise.
  Since $\hat\psi\in S(\R^d)$, we see that  the second integral in (\ref{j}) does not exceed
 $$
   c_N\int_{|\l-y|\ge1}(R|y-\l|)^{-N}|\hat\mu|(y)=c_N R^{-N}\int_1^\infty t^{-N}d|\hat\mu|(B(\l,t)).
   $$
 If $N$ is large enough, then $|\hat\mu|(B(\l,t))=O(t^{N-1})$ as $t\to\infty$. Arguing as in the previous lemma, we obtain that the second integral in (\ref{j}) tends to $0$ as $R\to\infty$.

 Furthermore, since the measure $\mu$ is translation bounded, we get for large $R$
 $$
  \left|\int_{\R^d}\psi(|x|/R)e^{-2\pi i\la\l,x\ra}d\mu(x)-\o_d^{-1}\int_{B(x,R)}e^{-2\pi i\la\l,x\ra}d\mu(x)\right|\le C\e R^d.
 $$
 Since $\e$ is arbitrary, we get the assertion of the Lemma.\bs

{\bf Remark}. We used  property of $\mu$ to be translation bounded only in the latter part of the proof.
\medskip

{\bf Proof of Theorem \ref{9}}. Suppose that $\mu$ is translation bounded measure with countable spectrum and slowly increasing measure $\hat\mu$.  By Lemma \ref{L1},  $\psi\star\mu(s)$ is almost periodic for every $\psi\in S(\R^d)$.
Check that $(\phi\star\mu)(t)$ is almost periodic for each continuous function $\phi$ with a compact support in the ball $B(R)$. Let  $\psi_n\in S(R^d),\ \supp\psi_n\subset B(R+1)$, be a sequence that uniformly converges to  $\phi$. Since $\mu$ is translation bounded,  the almost periodic functions $(\psi_n\star\mu)(t)$ uniformly converge to $(\phi\star\mu)(t)$, therefore the latter function is almost periodic, and $\mu$ is an almost periodic measure.

Now suppose that $\mu$ is almost periodic. By Theorem \ref{R}, $\mu$ is translation bounded. Let $\nu$ be the atomic component of measure $\hat\mu$, let $\mu_1$ be the inverse Fourier transform of $\nu$, and $\mu_2$ be the inverse Fourier transform of $\hat\mu-\nu$. It follows from the first part of the proof that $\mu_1$ is an almost periodic measure. Therefore, $\mu_2$ is almost periodic too. Let $\p$ be
$C^\infty$-function with compact support. The convolution $\p\star\mu_2$ is an almost periodic function, hence measure $\p\star\mu_2(x) m_d(x)$ is translation bounded. By Lemma \ref{L2}, Fourier coefficient of $\p\star\mu_2$ at each frequency $\l\in\R^d$ equals $(\p\star\mu_2)\hat\,(\{\l\})=\hat\p(\l)\hat\mu_2(\{\l\})=0$. Therefore, $\p\star\mu_2\equiv 0$ for each $\p\in S(\R^d)$ with compact support. Then $\mu_2=0$ and $\hat\mu=\nu$. \bs
\medskip

To prove Theorem \ref{3} we need two lemmas.
 \begin{Lem}\label{L3}
 Let $\mu=\sum_{\l\in\L}a_\l\d_\l,\ |a_\l|\ge\e>0$, be an almost periodic measure on $\R^d$ with uniformly discrete $\L$. Then the set $\L$ is almost periodic too.
\end{Lem}
 {\bf Proof}. Check that for every continuous function $\p(x)$ with compact support
 and every $\e>0$ the function $\sum_{\l\in\L}\p(x-\l)$ has a relatively dense set of $\e$-almost periods.
 Set
  $$
 \eta=(1/2)\inf_{\l,\l'\in\L,\l\neq\l'}|\l-\l'|,\quad  \b=\inf_{\l\in\L}|a_\l|.
  $$
     Without loss of generality suppose that $\supp\p\subset B(\eta)$, and  $0\le\p(x)\le \p(0)=1$. There exists $\rho\in (0,\eta)$ such that  the inequality $|x-x'|<\rho$ implies $|\p(x)-\p(x')|<\e$. Set $\psi(x)=\p(\eta x/\rho)$. Note that if $\psi(x)\neq0$, then $|x|<\rho$. Let $\tau$ be $\b$-almost period of the function $\psi\star\mu(x)=\sum_{\l\in\L} a_\l\psi(x-\l)$. We have for all $x\in\R^d$
   $$
       \left|\sum_{\l\in\L} a_\l\psi(x-\l)-\sum_{\l'\in\L} a_{\l'}\psi(x+\tau-\l')\right|<\b.
     $$
   Clearly, at most one term in every sum does not vanish. For $x=\l$ we get $|a_\l-a_{\l'}\psi(\l+\tau-\l')|<\b$, hence $\psi(\l+\tau-\l')\neq0$ and $|\l+\tau-\l'|<\rho$. Therefore, $|\p(x-\l)-\p(x+\tau-\l')|<\e$. If there is another $\l''\in\L$ such that
  $\p(x+\tau-\l'')\neq0$, then $|x+\tau-\l''|<\eta$ and $|\l'-\l''|<2\eta$. This is impossible, therefore for all $x\in\R^d$
  $$
       \left|\sum_{\l\in\L}\p(x-\l)-\sum_{\l'\in\L}\p(x+\tau-\l')\right|<\e.
  $$
Since there is a relatively dense set of $\b$-almost periods $\tau$ of the function $\psi\star\mu$, we conclude  that the function $\sum_{\l\in\L}\p(x-\l)$ is almost periodic too. \bs
\begin{Lem}\label{L4}
 Let $\mu=\sum_{\l\in\L}a_\l\d_\l$ be a uniformly discrete measure from $S'(\R^d)$ with the slowly increasing measure $\hat\mu$. Then $\sup_{\l\in\L}|a_\l|<\infty$ and, consequently, the measure $\mu$ is translating bounded.
\end{Lem}
{\bf Proof}. Set $\eta=(1/2)\inf_{\l,\l'\in\L,\l\neq\l'}|\l-\l'|$. Let $\psi(|y|)$ be a $C^\infty$-function  such that $\supp\psi(|y|)\subset B(\eta)$ and $\psi(0)=1$. We have
 $$
  \sup_{\l\in\L}|a_\l|\le\sup_{x\in\R^d}\left|\int\psi(|x-\l|)d\mu(\l)\right|=
  \sup_{x\in\R^d}\left|\int_{\R^d}\hat\psi(y)e^{2\pi i\la x,y\ra}d\hat\mu(y)\right|\le\int_{\R^d}|\hat\psi(y)|d|\hat\mu|(y).
 $$
 Taking into account  that $\hat\psi(y)\in S(R^d)$ and arguing as above in Lemma \ref{L1}, we obtain that the latter integral is
 finite. \bs
\medskip

Our proof of Theorem \ref{3} is based on the following result.
 \begin{Th}[\cite{F1}]\label{6}
  Let $\L$ be an almost periodic set in $\R^d$ with a discrete set of differences $\L-\L$. Then  $\L$ is a finite union of translates of a single full-rank lattice.
   \end{Th}
{\bf Proof of Theorem \ref{3}}. Since the set $\L-\L$ is discrete, we get that $\L$ is uniformly discrete. By Lemma \ref{L4}, $\mu$ is translation bounded. By Theorem \ref{9}, $\mu$ is almost periodic. Lemma \ref{L3} implies that $\L$ is almost periodic as well. Now, using Theorem \ref{6}, we obtain the assertion of Theorem \ref{3}. \bs
\medskip

There is an analog of Theorem \ref{6} for pairs of measures.
 \begin{Th}\label{8}
Let $\mu_1,\,\mu_2$ be large \F s such that the set $\L_1-\L_2$ of differences of their supports $\L_1,\,\L_2$ is discrete. Then there exists a full-rank lattice $L$ such that the both supports $\L_1,\,\L_2$ are finite unions of translates of $L$.
 \end{Th}
In fact, the proof of this theorem is the same as in the previous one. But instead of Theorem \ref{6} we need to use the corresponding assertion for pairs of measures from \cite{F2}.

\section{Wiener's Theorem for quasicrystals}

In this section we prove Theorems \ref{5} and \ref{7}. We begin with the following lemma.
\begin{Lem}\label{L5}
 Let $F(\t,\tau),\ \t=(\t_1,\dots,\t_N)\in [0,1]^N,\,\tau\in [0,1]$, be $C^\infty$-differentiable function in all variables and periodic with period $1$ in each $\t_j,\ j=1,\dots,N$. Then its Fourier series
  $$
  F(\t,\tau)=\sum_{n\in\Z^N}c_n(\tau) e^{2\pi i\la\t,n\ra}
  $$
 converges absolutely and uniformly in $\tau\in[0,1]$.
  \end{Lem}
{\bf Proof}. We have
$$
  c_n(\tau)=\int_{[0,1]^N}F(\t,\tau)e^{-2\pi i\la\t,n\ra}d\t.
$$
 Let $n=(n_1,\dots,n_N)$ with $n_1\neq0$. Integrating by parts in variable $\t_1$ two times, we get
$$
  c_n(\tau)=\frac{-1}{4\pi^2n_1^2}\int_{[0,1]^N}\frac{\partial^2F(\t,\tau)}{\partial\t_1^2}e^{-2\pi i\la\t,n\ra}d\t.
$$
If $n_1\cdot n_2\cdots n_m\neq0$, then after integrating by parts in variables $\t_1,\dots,\t_m$ we get
$$
  c_n(\tau)=\frac{(-1)^m}{(2\pi)^{2m}n_1^2\cdots n_m^2}\int_{[0,1]^N}\frac{\partial^{2m}F(\t,\tau)}{\partial\t_1^2\dots\partial\t_m^2}
  e^{-2\pi i\la\t,n\ra}d\t.
$$
 Therefore,
 $$
 |c_n(\tau)|\le\frac{1}{(2\pi)^{2m}n_1^2\cdots n_m^2}\max_{\t,\tau}\left|\frac{\partial^{2m}F(\t,\tau)}{\partial\t_1^2\dots\partial\t_m^2}\right|.
$$
Consequently, for all $n\in\Z^N$ we get
$$
|c_n(\tau)|\le C(F)\min\{1,n_1^{-2}\}\cdots\min\{1,n_N^{-2}\}.
$$
Obvious convergence of the series $\sum_{n\in\Z^k}\min\{1,n_1^{-2}\}\cdots\min\{1,n_N^{-2}\}$ implies the assertion of the Lemma. \bs

\medskip

{\bf Proof of Theorem \ref{5}}. Let $U$ be a neighborhood of $K$ such that the function $h(z)$ is holomorphic and bounded on $U$. Set $\eta=(1/4)\dist(K,\C\setminus U)$. Let $\p(|z|)$ be $C^\infty$-differentiable nonnegative function with support in $B(\eta)$ such that $\int_{B(\eta)}\p(|z|)dm_2(z)=1$. Then the function
$$
     H(z)=\int_{\R^d}h(z-\zeta)\p(|\zeta|)dm_2(\zeta)
$$
is $C^\infty$-differentiable on $\R^d$ and coincides with $h(z)$ on the set $\{z:\dist(z,K)\le3\eta\}$.

Suppose $f(x)=\sum_n c_ne^{2\pi i\la x,\g_n\ra}$. Choose $N<\infty$ such that for $S(x)=\sum_{n\le N}c_ne^{2\pi i\la x,\g_n\ra}$ we get $\|f(x)-S(x)\|_W<\eta$. Applying Lemma \ref{L5} to the function
$$
F(\t,\tau)=H\left(\sum_{n\le N}c_ne^{2\pi i\t_n}+2\eta e^{2\pi i\tau}\right)
$$
 and replacing $\t_j$ with $\la x,\g_j\ra,\ j=1,\dots,N$, we get
\begin{equation}\label{h1}
  H(S(x)+2\eta e^{2\pi i\tau})=\sum_n c_n(\tau)e^{2\pi i\la x,\rho_n\ra}\in W,\qquad \rho_n\in\Lin_\Z\{\l_1,\dots,\l_N\},
\end{equation}
and $\| H(S(x)+2\eta e^{2\pi i\tau})\|_W\le C(N)$ uniformly in $\tau\in[0,1]$.

Next, taking into account that $|f(x)-S(x)|\le\|f-S\|_W<\eta$, we have
\begin{equation}\label{h2}
(2\eta e^{2\pi i\tau}-[f(x)-S(x)])^{-1}=\sum_{k=0}^\infty\frac{[f(x)-S(x)]^k}{\left[2\eta e^{2\pi i\tau}\right]^{k+1}}.
\end{equation}
Clearly, this sum belongs to $W$ and its norm is uniformly bounded in $\tau\in[0,1]$. Hence the same assertion is valid for the
product of (\ref{h1}) and (\ref{h2}). Consequently, the function
$$
 g(x)=\int_0^1 \frac{H(S(x)+2\eta e^{2\pi i\tau})2\eta e^{2\pi i\tau}}{(2\eta e^{2\pi i\tau}-[f(x)-S(x)])}d\tau=
 \frac{1}{2\pi i}\int_{|\zeta-S(x)|=2\eta}\frac{H(\zeta)~d\zeta}{\zeta-f(x)}.
$$
 belongs to $W$.

 If $f(x)\in K$, then $\dist(S(x),K)<\eta$ and $B(S(x), 2\eta)\subset\{\zeta:\dist(\zeta,K)<3\eta\}$. Moreover, $f(x)\in B(S(x),2\eta)$ and the function $H(z)$ coincides with $h(z)$ for $z\in B(S(x),2\eta)$. Therefore, $g(x)=h(f(x))$ in this case.    \bs
\medskip

To deduce Theorem \ref{7} from Theorem \ref{5} we prove the following two lemmas.
\begin{Lem}\label{L7}
 Let $\mu$ be a slowly increasing measure with a translation bounded measure $\hat\mu$. Then for all $f\in W$ the Fourier transform $(f\mu)\hat\ $ is a  translation bounded measure.
\end{Lem}
{\bf Proof}. It is evident that $f\mu$ is a slowly increasing measure. Set  $\mu_\g(x)=e^{2\pi i\la x,\g\ra}\mu(x)$ for $\g\in\R^d$.  Clearly, for each $y\in\R^d$ we get $|\hat\mu_\g|(B(y,1))=|\hat\mu|(B(y-\g,1))$.
Suppose $f(x)=\sum_n c_ne^{2\pi i\la x,\g_n\ra}\in W$. Then $(f\mu)(x)=\sum_n c_n\mu_{\g_n}(x)$  and $(f\mu)\hat\ =\sum_n c_n\widehat\mu_{\g_n}(x)$. Therefore for each $y\in\R^d$
$$
|(f\mu)\hat\ |(B(y,1))\le \sum_n |c_n||\widehat\mu_{\g_n}|(B(y,1))\le \sum_n |c_n|\sup_{y\in\R^d}|\hat\mu|(B(y,1)).\quad \bs
$$
\begin{Lem}\label{L8}
Let $\mu=\sum_{\l\in\L}a_\l\d_\l$ be a uniformly discrete measure with countable spectrum and slowly increasing Fourier transform. Then there is a function $f\in W$ such that $f(\l)=a_\l$ for all $\l\in\L$.
\end{Lem}
{\bf Proof}. Set $\eta=(1/2)\inf_{\l,\l'\in\L,\l\neq\l'}|\l-\l'|$. Let $\psi$ be $C^\infty$-differentiable function with support in $B(\eta)$ such that $\psi(0)=1$. Then the function $f(x)=\psi\star\mu(x)=\sum_{\l\in\L}a_\l\psi(x-\l)$ satisfies the condition $f(\l)=a_\l$ for all $\l\in\L$. By Lemma \ref{L1}, $f\in W$. \bs

\medskip

{\bf Proof of Theorem \ref{7}}. We apply Theorem \ref{5} with $K=\overline{\{a_\l:\,\l\in\L\}}$ and $f\in W$ such that $f(\l)=a_\l$, which exists by Lemma \ref{L8}, and then we apply Lemma \ref{L7} for the measure $\mu$ and the function $g$ from Theorem \ref{5}.

\section{Generalization of Meyer's Theorem}

To prove the results of this section we need some  definitions.

 {\it Lattice} is a discrete subgroup of $R^d$. If $A$ is a lattice or a coset of some lattice in $\R^d$, then $\di A$ is the dimension of the smallest translated subspace of $\R^d$ that contains $A$.  Every lattice $L$ of dimension $k$ has the form $T\Z^k$, where $T:\,\Z^k\to\Z^d$ is a linear operator of rank $k$. For $k=d$ we get a  full-rank lattice.  Also, the coset ring of an abelian topological group $G$ is the smallest collection  of subsets of $G$, that is closed under finite unions, finite intersections and complements and which contains all cosets of all open subgroups of $G$.
Next, Bohr compactification $\hR$ of $\R^d$ is a compact group, its dual is $\R^d$ with the discrete topology,  $\R^d$ is a dense subset of $\hR$ with respect to the topology on $\hR$, and restrictions to $\R^d$ of continuous functions on $\hR$ are just  almost periodic functions on $\R^d$, in particular, they are bounded and continuous on $\R^d$ (see for example \cite{Ru}).

Both Meyer's Theorem \ref{1} and our Theorem \ref{4} are based on Cohen's Idempotent Theorem:
\begin{Th}[\cite{Co}]\label{Co}
Let $G$ be a locally compact abelian group and $\hat G$ its dual group. If $\nu$ is a finite Borel measure on $G$ and is such that
its Fourier transform $\hat\nu(\l)\in\{0,1\}$ for all $\l\in\hat G$, then the set $\{\l:\,\hat\nu(\l)=1\}$ is in the coset ring of $\hat G$.
\end{Th}
We will apply Theorem \ref{Co} to Bohr compactification $\hR$ of $\R^d$ and its  dual  $\R^d$ with the discrete topology. Then we will use the following theorem by M.Kolountzakis:
 \begin{Th}[\cite{K}] Elements of the ring of cosets of $\R^d$ in the discrete topology, which are discrete in the usual topology of $\R^d$, are precisely  finite unions of sets of the type
\begin{equation}\label{c}
A\setminus(\cup_{j=1}^N B_j),\quad A,\ B_j\ \mbox{discrete cosets},\quad\di B_j<\di A\quad \mbox{for all}\ j.
\end{equation}
\end{Th}

{\bf Proof of Theorem \ref{4}}.
Let $\rho$ be the inverse Fourier transform of the measure $\mu$, let $\psi(|x|)\in S(\R^d)$ be the function with compact support  from Lemma \ref{L2} with $\e=1/2$.

   Set $\rho_R(y)=R^{-d}\psi(|y/R|)\rho(y)$.  Integrating by parts and using  (\ref{m}), we get
    $$
   |\rho_R|\le R^{-d}\int_0^{(3/2)R}\psi(t/R)d|\rho|(B(t))\le R^{-d}\left[C'+C''R^{-1}\int_0^{(3/2)R} t^d|\psi'(t/R)|~dt\right].
  $$
  Therefore,
 $$
  \limsup_{R\to\infty}|\rho_R|\le C<\infty.
  $$
  Arguing as above in Lemma \ref{L2}, we get that $\hat\rho_R(x)$ tends to $\hat\rho(\{\l\})=\mu(\{\l\})=a_\l$ for $x=\l\in\L$ and to zero for
   $x\notin\L$.

 Since variations of the measures $\rho_R$ are uniformly bounded, they act on all bounded functions on $\R^d$, and hence also on all functions from $C(\hR)$. Therefore there exists a measure $\hr$ on $\hR$ with a finite total variation $|\hr|$, and a  subsequence $R'$ such that $\rho_{R'}\to\hr$ in the weak--star topology. In other words, $\la\rho_{R'},f\ra\to\la\hr,f\ra$ as $R'\to\infty$ for all $f\in C(\hR)$. Applying this to every character of $\hR$ in place of $f$ we obtain
 $$
 \hat\hr(x)=\lim_{R'\to\infty}\hat\rho_{R'}(x)=\begin{cases}a_\l, & x=\l\in\L,\\  0, & x\not\in\L.\end{cases}
 $$
 Note that $\hat\hr(x)$ is a continuous function with respect to the discrete topology on $\R^d$, and $|a_\l|\le|\hr|$ for all $a_\l$.

Using Lemma \ref{L8}, choose a function $f\in W$ such that $f(\l)=a_\l$ for all $\l\in\L$. Taking into account that $\inf_\l|a_\l|>0$ and using Theorem \ref{5} with $h(z)=1/z$, we  construct the function $g(x)=\sum_nc_ne^{2\pi i\la x,\tau_n\ra}\in W$ with the property $g(\l)=1/a_\l$. Since
$\hat\hr(y+\tau)=e^{2\pi i\la x,\tau\ra}\hat\hr$, we see that the Fourier transform of a finite sum $\sum_{n\le N}c_n\hr(y+\tau_n)$
is equal to $\sum_{n\le N}c_ne^{2\pi i\la x,\tau_n\ra}\hat\hr$.
 The variation of the measure $\sum_{n\le N}c_n\hr(y+\tau_n)$ does not exceed $\sum_{n\le N}|c_n||\hr|(\hR)$, therefore the previous assertion is valid for infinite sums. Hence the Fourier transform of the measure $\frak b=\sum_nc_n\hr(y+\tau_n)$ equals $1$ for $x=\l\in\L$ and $0$  for $x\not\in\L$.
By Theorem \ref{Co}, $\L$ is in the coset ring of $\R^d$ in discrete topology. Taking into account Lemma  \ref{L4} and Theorem \ref{9}, we see that the measure $\mu$ is almost periodic. Therefore, by Lemma \ref{L3},  measure $\sum_{\l\in\L}\d_\l$ is also almost periodic.

Using (\ref{c}) and the notion of almost periodicity, it is easy to get the assertion of our theorem.  By definition, put $\d_A=\sum_{x\in A}\d_x$ for a countable set $A\subset\R^d$. We have
$$
\L=\cup_{k=1}^K A_k\cup(\cup_{j=1}^J B_j)\setminus\cup_{i=1}^I C_i,
$$
where $A_k,\,B_j,\,C_i$ are cosets such that $\dim A_k=d,\,\dim B_j<d,\,\dim C_i<d$. Let $A_1=x_1+L_1,\ A_2=x_2+L_2$ be two cosets of dimension $d$. If $\dim A_1\cap A_2=d$, then $A_1\cap(\R^d\setminus A_2)$ is a finite union of disjoint cosets of $L_1\cap L_2$, and the same is valid for $A_2\cap(\R^d\setminus A_1)$. Therefore we can replace $A_1\cup A_2$ by a finite sum of disjoint cosets $A'_l,\,\dim A'_l=d$, hence $\d_{A_1\cup A_2}=\sum_l\d_{A'_l}$. For every two cosets $H_1,\,H_2$ with $H_1\cap H_2\neq\emptyset$ and $\dim(H_1\cap H_2)<d$ we get $\d_{H_1\cup H_2}=\d_{H_1}+\d_{H_2}-\d_{H_1\cap H_2}$.  Repeating this transformation for each pair of cosets with nonempty intersection, we  obtain a representation of $\d_\L$ after a finite number of steps of the form
$$
   \d_\L=\sum_{s=1}^{N_1}\d_{D_s}+\sum_{s=1}^{N_2}\d_{E_s}-\sum_{s=1}^{N_3}\d_{F_s},
$$
where  $D_s,\,E_s,\,F_s$ cosets, $\dim D_s=d$, $\dim E_s<d$, and $\dim F_s<d$. Obviously, every coset $D_s$ has $d$ linearly independent periods, hence $\d_{D_s}$ is an almost periodic measure, and so is $\sum_{s=1}^{N_1}\d_{D_s}$. Therefore the measure
 \begin{equation}\label{b}
   \sum_{s=1}^{N_2}\d_{E_s}-\sum_{s=1}^{N_3}\d_{F_s}=\d_\L-\sum_{s=1}^{N_1}\d_{D_s}
 \end{equation}
is almost periodic too. But its support is contained in a finite union of hyperplanes and isn't relatively dense. Therefore the measure (\ref{b}) is identically zero, and $\d_\L=\sum_{s=1}^{N_1}\d_{D_s}$. Since  the measure $\d_\L$ has the unit masses at every point of $\L$, we see that the cosets $D_s$ are disjoint and $\L=\cup_{s=1}^{N_1}D_s$.  \bs
\medskip

{\bf Proof of the Corollary}. By Lemma \ref{L4}, $A$ is a bounded set.  Apply Theorem \ref{7} with a function $h(z)$  that equals $1$ on a neighborhood $U$ of $A$ and $0$ on $\C\setminus U$, and then Theorem \ref{4} for $\sum_{\l\in\L, a_\l\in A}a_\l\d_\l$. \bs

\end{document}